\documentclass{article}

\usepackage{amsmath,amsfonts,amssymb}
\usepackage{mathptmx}       
\usepackage{helvet}         
\usepackage{courier}        
\usepackage{a4wide}        
 \usepackage{authblk}
\usepackage{color}
\usepackage{todonotes}
\usepackage{makeidx}         
\usepackage{graphicx}        
\usepackage{multicol}        
\usepackage[bottom]{footmisc}

\newcommand{\compliance}{\kappa}
\newcommand{\jmp}[1]{[\![#1]\!]}
\newcommand{\av}[1]{\{\!\!\!\{#1\}\!\!\!\}_w}
\newcommand{\avd}[1]{\left<\!\!\left<#1 \right>\!\!\right>_w}
\newcommand{\lom}[2]{(#1,#2)_\Omega}
\newcommand{\lga}[2]{\left<#1,#2\right>_\Gamma}

\date{}
\begin{document}

\title{Deriving robust unfitted finite element methods from augmented Lagrangian formulations}
\author[$\dagger$]{Erik Burman}
\author[$\star$]{Peter Hansbo}
\affil[$\dagger$]{\small Department of Mathematics, University College London, London, UK--WC1E  6BT, United Kingdom}
\affil[$\star$]{\small Department of Mechanical Engineering, J\"onk\"oping University, SE-55111 J\"onk\"oping, Sweden}

%
%
\maketitle

\abstract{In this paper we will discuss different coupling methods {suitable for use in} the framework of the recently introduced
  CutFEM paradigm, cf. Burman et al. \cite{BuClHaLaMa15}.
In particular we will consider mortaring using Lagrange multipliers on
the one hand and Nitsche's method on the other. For simplicity we will first discuss these method in the setting of
uncut meshes, and end with some comments on the extension to CutFEM.
We will, for comparison, discuss some different types of problems such as high contrast
problems and problems with stiff coupling or adhesive contact. We will review some of
the existing methods for these problems and propose some alternative
methods resulting from crossovers from the Lagrange multiplier framework to
Nitsche's method and vice versa.}

\section{Introduction}
\label{sec:1}
Recently there has been increased interest in unfitted finite element
for the imposition of boundary conditions or more generally for the coupling of
physical systems over an interface. The unfitted discretization is then made as
independent as possible of the geometric description of interfaces and domain boundaries
in order to minimize the
complexity of mesh generation. 
One such method is the Cut Finite Element Method (CutFEM) \cite{BuClHaLaMa15},
the goal of which is to retain the accuracy and robustness of a
standard finite element method. To reach this aim, stabilization techniques are applied
to make both the accuracy of the approximation and the system
condition number independent of the mesh/boundary intersection and
physical parameters. Thanks to this robustness of the discretization, 
powerful linear algebra techniques developed for finite element
methods are made available for solving the linear systems
obtained by the CutFEM discretization. 

In the CutFEM approach the boundary of---or interfaces in---a given domain is represented on a background 
grid, e.g., using a level set function, and the background grid is 
also used to discretize the governing partial differential 
equations. CutFEM builds on a general finite element formulation for the approximation 
of partial differential equations, in the bulk and on surfaces (interfaces or boundaries), that can handle elements 
of complex shape and where boundary and interface conditions are built
into the discrete formulation. CutFEM requires only
a low-quality, even non-conforming, surface mesh representations of the
computational geometry, thus mitigating mesh generation work.

Unfitted methods typically use either
Lagrange multipliers or Nitsche's method for the mortaring on interfaces or boundaries and it has
been shown that Nitsche's method can be derived from a stabilized
Lagrange muliplier method, due to Barbosa and Hughes \cite{BaHu91} after
static condensation of the multiplier \cite{St95, Ju15}, cf. Section \ref{Nitstab}. 

Another,
fluctuation based, stabilization of the Lagrange multiplier was
proposed in \cite{BuHa10a} and further developed in the works
\cite{BuHa10,Bu14,BC12} and also \cite{LF17} of this collection. For this type of methods the multiplier
typically may not be eliminated, which can be a disadvantage; however,
in some situations it is desirable to define the multiplier on a
different space anyway. 

Regardless of the relative virtues of the two
approaches it is interesting to compare them and see how developments
for one type of methods can be exploited in the context of the
other. Similarly as in \cite{St95} we will here be interested in
deriving methods using Lagrange multipliers and then recover the
associated Nitsche method by formal elimination of the
multiplier. However in our case we will base the discussion on the
concept of augmented Lagrangian methods, {which has} recently  been
successfully applied in the context of contact problems using
Nitsche's method \cite{ChHi13,ChHiRe15} and Lagrange multipliers
\cite{BuHaLa16}. 

The idea behind the augmented Lagrangian is to add a
least squares term on the constraint to the Lagrangian functional of
the constrained optimization problem on the discrete level. Provided this modification is
not too strong it is expected to improve the conditioning as well as
improving the control of the constraint compared to the case where
only the multiplier is used to enforce the constraint. For early work
on augmented Lagrangian methods in computational methods for partial
differential equations we refer to Glowinski and le Tallec \cite{GT82}
or Fortin and Glowinski \cite{FG83}. However, in
cases where the problem depends on physical parameters that have some
singular behavior, the augmented Lagrangian can lead to a severely
ill-conditioned problem. This is typically the case for problems with
high contrast in the diffusivity or strongly bonded adhesive
problems. The ill-conditioning results from the fact that in the singular
limit the least squares term on the constraint blows up, leading to
ill-conditioning and also locking, unless the mesh-size is small enough
to resolve the singularity. In many problems such a resolution is
unfeasible and it is then useful to relax the size of the least squares 
term. This idea has been implicitly used in a number of works,
starting with the paper \cite{HaHa04} on debonding problems
using Nitsche's method on unfitted meshes and then further developed
for free flow porous media coupling in \cite{BuHa07} and boundary
conditions in \cite{JuSt09}.  A recent inventive use of this framework 
was proposed in \cite{FL15} where it was  applied to coupling schemes in fluid-structure interaction. 
All of these works are concerned with Nitsche type
formulations.

Nevertheless there seems to have been no {attempts at exploring}
these ideas directly in the framework of augmented
Lagrangian methods. Our main objective in this paper is to study some
model problems, exhibiting the typical parameter dependent behaviour,
in the framework of augmented Lagrangian methods.

The idea is to first revisit the imposition of boundary conditions
with Lagrange multipliers and Nitsche's method.
Then we consider three different situations of
domain decomposition. Here we assume that the problem
is set on a domain $\Omega$ with two systems, defined in the subdomains
$\Omega_1$ and $\Omega_2$, that are coupled over a
smooth boundary $\Gamma$. 
For simplicity we consider Poisson type
problems and do not discretize the boundary $\Gamma$, that is, we
consider the semi-discretized setting. The discussion can easily be
extended to for instance compressible elasticity.
We will consider the following three model problems:
\begin{enumerate}
\item general boundary conditions;
\item Poisson's equation with high contrast in the diffusion coefficient;
\item debonding and adhesive contact.
\end{enumerate}
First, in Section \ref{sec:Nitsche}, we will recall, for an uncut mesh with Dirichlet boundary conditions, how least squares stabilized Lagrange multiplier methods lead to 
Nitsche's method following \cite{St95}, and show how, alternatively, the augmented Lagrangian approach leads to the
same formulation. Then, in Sections \ref{sec:gen_bc}--\ref{sec:adhesive} we consider the three
different model problems in the augmented Lagrangian framework and
derive robust Nitsche methods as well as robust Lagrange multiplier
methods, still formulated on uncut meshes.
Some of the proposed methods appear to be new, whereas others are
known in the literature and we will discuss existing results for the
methods and without going into technical details we will
speculate on what results are likely to carry over to the cases
considered herein. In Section \ref{sec:cutmesh}, we round off with some remarks considering the extension to cut finite element meshes, in particular with respect to
stabilization of the discrete system, and, finally, in Section \ref{sec:example}, we give a numerical example for one of the model problems.

\section{Derivations of Nitsche's method from Lagrange multipliers\label{sec:Nitsche}}
\subsection{Model problem}
Let us first consider the typical Poisson model problem
of finding $u$ such that
\begin{equation}\label{poiss}
-\Delta u= f ~\text{in} ~\Omega ,\quad u=g ~ \text{on} ~\Gamma:=\partial\Omega,
\end{equation}
where $\Omega$ is a bounded domain in two or three space
dimensions, with outward pointing normal $\boldsymbol{n}$,
and $f$ and $g$ are given functions. For simplicity, we shall assume that $\Omega$ is polyhedral (polygonal).  
The typical way of prescribing $u=g$ on the boundary is to pose the problem (\ref{poiss}) as a minimization
problem with side conditions and seek stationary points to the functional
\begin{equation}\label{eq:first}
\mathcal{L}(v,\mu) := \frac12 a(v,v)  -
\lga{\mu}{v-g}  - \lom{f}{v} ,
\end{equation}
where 
\[
\lom{f}{v} := \int_{\Omega}f v \, d\Omega, \; a(u,v) := \int_\Omega \nabla u\cdot\nabla v\, d\Omega, \; \lga{\mu}{v-g} := \int_{\Gamma}\mu (v-g) \, ds .
\]
The stationary points are given by finding $(u,\lambda)\in H^1(\Omega)\times H^{-1/2}(\Gamma)$ such that
\begin{equation}\label{eq:weakform}
a(u,v)-\lga{\lambda}{v} = (f,v)\quad \forall v\in H^1(\Omega),
\end{equation}
\begin{equation}\label{eq:weakform2}
\lga{\mu}{u} = \lga{\mu}{g}\quad \forall \mu\in H^{-1/2}(\Gamma) .
\end{equation}
As is well known, the discretization of this problem requires
balancing of the discrete spaces for the multiplier $\lambda$ and the primal solution $u$ in order for the method to be stable, for examples cf. \cite{BrFo91}. In the following we shall not consider balanced ({\em inf-sup}\/ stable) discrete methods but instead focus on stabilized methods.

\subsection{Nitsche's method as a stabilized multiplier method\label{Nitstab}}
Formally, the Lagrange multiplier in (\ref{eq:weakform}) is given by $\lambda = \partial_n u$, where $\partial_n{v}:={\boldsymbol n}\cdot\nabla v$, and a well known stabilization method \cite{BaHu91} for the discretization of (\ref{eq:weakform})--(\ref{eq:weakform2}) to is to add
a term penalizing the difference between the discrete multiplier and the discrete normal derivative of the primal solution. 
To this end, we assume that $\mathcal{T}_h$ is a
conforming shape regular meshe on $\Omega$, consisting of
triangles $T$
and 
define $V_h$ as the space of $H^1$--conforming piecewise polynomial
functions 
on $\mathcal{T}$,
\[
V_h := \{v_h \in H^1(\Omega): v_h\vert_T \in \mathbb{P}_k(T), \,
\forall T \in \mathcal{T} \},\quad \mbox{ for } k \ge 1.
\]
As discrete space for the multiplier, we define the trace mesh on $\Gamma$ as the set $\mathcal{F}_h$ of element faces $F$ on $\Gamma$
and set
\[
\Lambda_h := \{q_h \in L_2(\Gamma): q_h\vert_F \in \mathbb{P}_l(F), \,
\forall F \in \mathcal{F} \},\quad \mbox{ for } l \ge 0.
\]
Then 
we seek $(u_h,\lambda_h)\in V_h\times\Lambda_h$ such that
\begin{equation}\label{discreteBH}
a(u_h,v)-\lga{\lambda_h}{v}-\lga{\mu}{u_h}-\frac{1}{\gamma_0}{\lga{h(\lambda_h-\partial_n u_h)}{\mu-\partial_n v}} = (f,v)-\lga{\mu}{g}
\end{equation}
for all $(v,\mu)\in V_h\times\Lambda_h$. Here $h$ is the meshsize of the trace mesh on $\Gamma$, interpreted as a piecewise constant function along $\Gamma$ and $\gamma_0$ is a number to be chosen sufficiently large to obtain a stable method.
Following Stenberg \cite{St95} we now let $P_h:L_2(\Gamma)\rightarrow \Lambda_h$ denote the $L_2$--projection, and considering $\Lambda_h$ to be a space of discontinuous
discrete multipliers, continuous inside each face $F$ of the Lagrange multiplier mesh on $\Gamma$, we can eliminate the multiplier from (\ref{discreteBH}):
\begin{equation}
\lambda_h\vert_E = P_h\partial_n u_h\vert_E -\gamma_0 h^{-1}P_h (u_h-g)\vert_E\quad \forall E.
\end{equation}
Now considering the limiting case of $\Lambda_h \rightarrow L_2(\Gamma)$ we see that
\[
\lambda_h \rightarrow \partial_n u_h -\gamma_0 h^{-1} (u_h-g)
\]
and we can reintroduce this multiplier into our stabilized method,
{replacing also $\mu$ by $\partial_n v - \gamma_0
  h^{-1} v_h$,} to obtain the problem of finding $u_h\in V_h$ such that
\begin{equation}\label{Nit1}
a(u_h,v)-\lga{\partial_n u_h}{v}-\lga{\partial_n v_h}{u_h}+\gamma_0\lga{h^{-1} u_h}{v} = L(v)\quad\forall v\in V_h,
\end{equation}
where
\[
L(v) :=(f,v)-\lga{\partial_nv}{g}+\gamma_0\lga{ h^{-1}v}{g}
\]
which is Nitsche's method \cite{Nit70}.

\subsection{Nitsche's method as an augmented Lagrangian method}
The other approach to deriving Nitsche's method from Lagrange multipliers is more in the vein of Nitsche' original paper \cite{Nit70} where the method was
derived from a discrete minimization problem without multipliers. The 
Lagrangian in (\ref{eq:first}) on the discrete spaces is augmented by a penalty term mutliplied by $\gamma\in \Bbb{R}^+$ so that we seek stationary points to
\begin{equation}\label{eq:aug1}
\mathcal{L}(v,\mu) := \frac12 a(v,v)  -
\lga{\mu}{v-g}  +\frac12\|\gamma^{1/2} (v-g)\|^2_\Gamma - \lom{f}{v} ,
\end{equation}
leading to the problem of finding $(u,\lambda)\in H^1(\Omega)\times H^{-1/2}(\Gamma)$ such that
\[
 a(u,v)  -
\lga{\lambda}{v}  +\lga{\gamma\, u}{v} -\lga{\mu}{u} =  \lom{f}{v}+\lga{\gamma\, g}{v}-\lga{\mu}{g}  \\
\]
for all $(v,\mu)\in H^1(\Omega)\times H^{-1/2}(\Gamma)$. Choosing now in the discrete case $\lambda_h := \partial_n u_h$, $\mu = \partial_n v$, and $\gamma = \gamma_0h^{-1}$
we recover (\ref{Nit1}). It should be noted that augmented Lagrangian methods are not in general {\em inf-sup}\/ stable; typically an unstable method remains unstable and the augmentation rather serves the purpose of strengthening a method where the side condition is too weakly enforced. It is therefore rather remarkable that
in the particular case where the
discrete multiplier is replaced by the discrete normal derivative, the
augmentation always works as a stabilization mechanism, see also the
discussion in \cite{Bu14}.

An important feature of the augmented Lagrangian approach is that it directly carries over to the
case of inequality constraints, as first shown by Chouly and Hild in the context of elastic contact \cite{ChHi13}. In our model problem
we replace the constraint $u=g$ on $\Gamma$ by an inequality constraint $u -g \leq 0$ on $\Gamma$. We then have the following Kuhn--Tucker conditions on the multiplier and side condition:
\begin{equation}\label{eq:Kuhn1}
u -g \leq 0,\quad \lambda \leq 0, \quad \lambda (u-g)=0.
\end{equation}
The key to incorporating these conditions into the augmented Lagrangian scheme, as pioneered by Alart and Curnier \cite{AC91}, is to make the observation that
(\ref{eq:Kuhn1}) is equivalent to 
\begin{equation}\label{eq:lambda1}
\lambda = -{\gamma}\,[u-g-\gamma^{-1}\, \lambda]_+ \mbox{ and } u - g = [u-g-\gamma^{-1}\, \lambda]_-
\end{equation}
where $\gamma \in \mathbb{R}^+$, $[x]_{\pm}=\pm\max (\pm x,0)$, cf. \cite{ChHi13}.  Denoting
$P_\gamma(u,\lambda) := \gamma\,  (u-g)-\lambda$ we see that 
\begin{equation}\label{eq:pivot}
\lambda = -[P_\gamma(u,\lambda)]_+ =
[P_\gamma(u,\lambda)]_- -  P_\gamma(u,\lambda).
\end{equation}
We then formally write the
augmented Lagrangian, similar to \eqref{eq:aug1}, but
  using the second relation in (\ref{eq:lambda1}) for the contact constraint,
\begin{multline}\label{eq:aug2}
\mathcal{L}(v,\mu) := \frac12 a(v,v)  -
\lga{\mu}{v-g- [v-g-\gamma\, \mu]_-}  \\
+\frac12 \gamma \|v-g- [v-g-\gamma^{-1}\, \lambda]_-\|^2_\Gamma - \lom{f}{v}.
\end{multline}
Observe that this naive formulation is not differentiable, so of
little practical
use. Our aim is now to propose a modified Lagrangian that is
consistent with \eqref{eq:aug2}, but smoother. We use the relations
\eqref{eq:Kuhn1} for the modifications, in order not to perturb the
stationary points. The modifications must also respect the saddle
point structure of the system. First observe that by the two relations
of \eqref{eq:lambda1} the stationary point must satisfy $\lga{\mu}{
  [v-g-\gamma\, \mu]_-} = 0$, using this in the second term in the
right hand side of \eqref{eq:aug2} and developing the square of the
third term we obtain
\begin{multline*}
\mathcal{L}(v,\mu) := \frac12 a(v,v)  -
\lga{\mu}{v-g}   
+\frac12 \gamma \|v-g\|^2\\ + \frac12 \gamma^{-1} \| [P_\gamma(u,\mu) ]_-\|^2_\Gamma -\lga{(v-g)}{ [P_\gamma(u,\mu )]_-}- \lom{f}{v}.
\end{multline*}
We then use that the stationary point must satisfy the second relation
of \eqref{eq:lambda1} implying that
\[
\lga{(v-g)}{ [P_\gamma(u,\mu )]_-} = \gamma^{-1} \| [P_\gamma(u,\mu) ]_-\|^2_\Gamma
\]
resulting in 
\begin{multline*}
\mathcal{L}(v,\mu) := \frac12 a(v,v)  -
\lga{\mu}{v-g}   
+\frac12 \gamma \|v-g\|^2\\ - \frac12 \gamma^{-1} \| [P_\gamma(u,\mu) ]_-\|^2_\Gamma - \lom{f}{v}.
\end{multline*}
This formulation is differentiable and the corresponding Euler-Lagrange
equations read
\begin{multline*}
 a(u,v)  -
\lga{\lambda}{v}  +\lga{\gamma\, u}{v} -\lga{\mu}{u} \\
-\gamma^{-1} \lga{[\gamma (u-g)-\partial_n u]_-}{\gamma\,  v-\mu}=  \lom{f}{v}+\lga{\gamma\, g}{v}-\lga{\mu}{g}.  
\end{multline*}
Choosing now in the discrete case {$\lambda_h := \partial_n u_h$ and $\mu = \partial_n v$}
results in a formulation that
equals \eqref{Nit1} up to a nonlinear perturbation: we seek $u_h\in V_h$ such that
\begin{multline}\label{Nitcontact}
a(u_h,v)-\lga{\partial_n u_h}{v}-\lga{\partial_n
  v_h}{u_h}+\lga{\gamma u_h}{v} \\
-{\gamma}^{-1} \lga{\,[\gamma (u_h-g)-\, \partial_n u_h]_-}{\gamma\,  v-\partial_n v} = L(v)\quad\forall v\in V_h,
\end{multline}
where
\[
L(v) :=(f,v)-\lga{\partial_nv}{g}+{\gamma\lga{v}{g}}
\]
for all $v\in  V_h$.

To see the equivalence of this formulation with that introduced by
Chouly and Hild \cite{ChHi13} we once again use the relation
\eqref{eq:pivot} to write
{\begin{multline*}
-{\gamma}^{-1} \lga{\,[\gamma (u_h-g)-\, \partial_n u_h]_-}{\gamma\,  v-\partial_n v} = \\
-\lga{
 [P_\gamma(u,\partial_n u) ]_- -  P_\gamma(u,\partial_n u)
  }{ v-\gamma^{-1} \partial_n v}  -  \lga{ P_\gamma(u,\partial_n u)
  }{ v-\gamma^{-1}\partial_n v}\\
= \lga{
 [P_\gamma(u,\partial_n u) ]_+}{ v-\gamma^{-1} \partial_n v}-  \lga{ P_\gamma(u,\partial_n u)
  }{ v-\gamma^{-1}\partial_n v}.
\end{multline*}
}
Applying this relation in the last term of the left hand side of
\eqref{Nitcontact} and simplifying results in 
we seek $u_h\in V_h$ such that
\begin{equation}\label{eq:contact_Nitsche}
a(u_h,v)+\lga{{\gamma}\,[u_h-g-\gamma^{-1}\, \partial_n u_h]_+}{v-\gamma^{-1}\, \partial_n v} - \lga{\gamma^{-1}\,\partial_n u_h}{\partial_n v}= (f,v),
\end{equation}
for all $v\in  V_h$. With the choice $\gamma = \gamma_0/h$ (\ref{eq:contact_Nitsche})
precisely coincides with the nonlinear Nitsche method proposed by
Chouly and Hild \cite{ChHi13}. This variational problem corresponds to stationarity of the Lagran\-gian
\[
\mathcal{L}(v,\mu) := \frac12 a(v,v)  +
\frac12 \Vert \gamma^{1/2}[v-g-\gamma^{-1} \mu]_+\Vert_\Gamma^2 -\lga{\gamma^{-1}\mu}{\mu} - \lom{f}{v} ,
\]
cf. \cite{AC91}.

In the following we shall explore how the augmented Lagrangian approach can be used in different situations.

\section{Poisson's equation with general boundary conditions}\label{sec:gen_bc}
We first consider the simple case of imposing general boundary
conditions in a finite element method with the mesh fitted to the
boundary. Following \cite{JuSt09} we consider the model problem
\begin{align}
-\nabla \cdot \left(\varepsilon \nabla u\right) & = f \mbox{ in } \Omega\\
\varepsilon \partial_n u & = \compliance^{-1} (u_0 - u)+g \mbox{ on } \Gamma=\partial \Omega ,
\end{align}
where $\varepsilon$, $u_0$, and $\compliance$ are given functions of sufficient regularity.
The problem can be equivalently formulated as the stationary points of
the Lagrangian
\[
\mathcal{L}(v,\mu) := \frac12 a(v,v)  -
\lga{\mu}{v}  - \frac12 \compliance
\|\mu^2\|_{\Gamma}- \lom{f}{v} +\lga{(u_0 + \compliance g)}{\mu}.
\]
where 
\[
a(u,v) := (\varepsilon \nabla u, \nabla v)_\Omega.
\]
The augmented Lagrangian is obtained by adding a scaled least squares
term on the constraint. This leads to
\begin{align}\nonumber
\mathcal{L}(v,\mu) := {}&\frac12 a(v,v) -
\lga{\mu}{v}  - \frac12 \compliance
\|\mu\|^2_{\Gamma} \\ &
+ \frac{1}{2\compliance}\|\compliance \mu+ v -
u_0 - \compliance g\|_\Gamma^2- \lom{f}{v} +\lga{(u_0 + \compliance g)}{\mu}.
\end{align}
Developing the square of the least squares contribution leads to the
elimination of the multiplier and the form
\begin{equation}
\mathcal{L}(v) := a(v,v) 
 - \frac{1}{2 \compliance}
\|v\|^2_{\Gamma} - \lom{f}{v} - \lga{(\compliance^{-1}u_0 + g)}{v}.
\end{equation}
The corresponding optimality system is the standard weak formulation
of the Robin problem: find $u \in H^1(\Omega)$ such that
\begin{align}
a(u,v) + \compliance^{-1} \lga{u}{v} = \lom{f}{v} + \lga{(\compliance^{-1} u_0 + g)}{v}, \quad \forall v
  \in H^1(\Omega).
\end{align}
Restricting the formulation to finite dimensional spaces results in 
the finite element formulation: find $u_h \in V_h$
such that
\begin{align}
a(u_h,v_h) + \compliance^{-1} \lga{u_h}{v_h} = \lom{f}{v_h} + \lga{(\compliance^{-1} u_0 + g)}{v_h}.
\end{align}
We see that as the ``compliance'' $\compliance \rightarrow 0 $, $u\vert_{\Gamma} \rightarrow
u_0$, implying a Dirichlet condition in the limit and as $\compliance \rightarrow
\infty$ we recover the weak formulation for the Neumann problem with
$\varepsilon \partial_n u\vert_\Gamma = g$.
As was pointed out in \cite{JuSt09} the corresponding linear system
becomes ill-posed if $\compliance<< h$. As a remedy for this we will temper the
coefficient in front of the least squares term in the augmented Lagrangian. Indeed if we introduce 
\[
S_h = (\compliance + h/{\gamma_{\compliance}})^{-1},
\]
where ${\gamma_{\compliance}}$ is a free parameter, and use this as coefficient, we
obtain the augmented Lagrangian
\begin{align}\nonumber
\mathcal{L}(v,\mu) := {}&\frac12 a(v,v) -
\lga{\mu}{v}  - \frac12 \compliance
\|\mu\|^2_{\Gamma} \\
& + \frac{S_h}{2}\|\compliance \mu+ v - u_0 -
\compliance g\|_\Gamma^2- \lom{f}{v} - \lga{(u_0 + \compliance g)}{\mu}.
\end{align}
The cancellations that led to the previous simple optimality system
will not take place this time as the moderated parameter $S_h$ will
prevent the system from becoming too stiff. Deriving we find the
following optimality system
\begin{multline}\label{eq:lagrange_mult}
a(u,v) - \lga{(1- \compliance S_h) \lambda}{v}-\lga{(1- \compliance S_h)
  \mu}{u}\\ - \lga{\compliance(1-S_h \compliance) \lambda}{\mu}
+\lga{S_h u}{v} = \\ \lom{f}{v}+ \lga{(u_0 + \compliance g)}{(1- S_h \compliance)\mu+S_h v}.
\end{multline}
Using this formulation in a finite element method with a pair $V_h
\times \Lambda_h$ that satisfies the {\em inf-sup}\/ condition leads to a
robust and accurate method. 
Formally replacing $\lambda$ by $\varepsilon \partial_n u$ and
$\mu$ by $\varepsilon \partial_n v$
and writing the form on the finite space $V_h$ yields the method: find $u_h \in
V_h$ such that
\begin{multline}
a(u_h,v_h) - \lga{(1- \compliance S_h) \varepsilon \partial_n u_h}{v_h}-\lga{(1- \compliance S_h)
 \varepsilon \partial_n v_h}{u_h}\\ - \lga{\compliance(1-S_h \compliance) \varepsilon \partial_n u_h}{\varepsilon \partial_n v_h}
+\lga{S_hu_h}{v_h} \\
= \lom{f}{v_h}+ \lga{(u_0 +\compliance g)}{(1-S_h \compliance)\varepsilon \partial_n v_h+S_h v_h}.
\end{multline}
We identify the Nitsche method proposed in \cite{JuSt09} and conclude that their
analysis is valid for $u_h$ if the parameter ${\gamma_{\compliance}}$ is chosen large
enough.
Observe that the Lagrange multiplier formulation
\eqref{eq:lagrange_mult}, which appears to be new, also varies
robustly between a Neumann and a Dirichlet condition in the two
limits, without succumbing to the ill-conditioning in the limits of
high or low  $\compliance$. We note in passing that sometimes the multiplier
method can have advantages compared to Nitsche's method, in particular if
fields defined on different meshes must be coupled, see for
instance \cite{BuHa10a}.
\section{FEM for Poisson's problem with discontinuous coefficients}
\label{sec:2}
We are interested in the following problem: find $u_i:\Omega_i \mapsto
\mathbb{R}$, $i=1,2$, such that
\begin{align}
-\nabla \cdot \varepsilon_i \nabla u_i & = f \mbox{ in } \Omega_i, \; i=1,2\\
u_i&=0 \mbox{ on } \partial \Omega\cap\Omega_i \; i=1,2\\
\jmp{\varepsilon \partial_n u}&=g\mbox{ on } \Gamma :=\partial\overline{\Omega}_1\cap\partial\overline{\Omega}_2\\
\jmp{u}&=0 \mbox{ on } \Gamma.
\end{align}
We use the notation $u=(u_1,u_2)\in V_1\times
V_2$ with the continuous spaces
\[
V_{i}=\left\{v_i\in H^1(\Omega_{i}): \
\partial v_i/\partial n_i  \in L_2(\Gamma), \
 v_i\vert_{\partial\Omega\cap\partial\Omega_{i}}=0\right\} ,\quad
i=1,2.
\]
Then $\jmp{u}$ denotes the jump of $u$ over $\Gamma$
defined as
\[
\jmp{u} := \lim_{\epsilon \rightarrow 0^+} u(x-\epsilon n) - u(x+\epsilon n)
\]
for $x \in \Gamma$ and $n$ denoting the normal on $\Gamma$ pointing
from $\Omega_1$ to $\Omega_2$. The diffusion coefficients $\varepsilon_i$, $i=1,2$, are assumed to be
constant functions. We will also use the weighted averages
\[
\av{u} := \lim_{\epsilon \rightarrow 0^+} ( w_1 u(x-\epsilon n)+
  w_2 u(x+\epsilon n))
\]
and
\[
\avd{u} := \lim_{\epsilon \rightarrow 0^+} ( w_2 u(x-\epsilon n)+
  w_1 u(x+\epsilon n))
\]
where $w_1,\, w_2 \ge 0$ are positive weights such that $w_1 + w_2 = 0$.
This problem can be shown to be equivalent to finding $(u,\lambda)
\in H^1(\Omega_1\cup\Omega_2) \times H^{-\frac12}(\Gamma)$, the
saddle point of the constrained minimization problem defined by the Lagrangian
\[
\mathcal{L}(v,\mu) := \frac12 a(v,v) +
\left<\mu,\jmp{v} \right>_{-\frac12,\frac12,\Gamma}  - \lom{f}{v} -\lga{g}{\avd{v}}.
\]
where now
\[
a(u,v) := (\varepsilon \nabla u,\nabla v)_{\Omega_1 \cup \Omega_2}.
\]
Thus, $(u,\lambda)$ fulfills
\begin{equation}\label{infsupref}
\mathcal{L}(u,\lambda)  =  \inf_{v \in V} \sup_{\mu \in \Lambda} \mathcal{L}(v,\mu).
\end{equation}
We also know that for the exact solution there holds $\lambda = -\av{\varepsilon
\nabla u \cdot n}$, for any admissible weights $w_1,w_2$.

To formulate a discrete method, we suppose that we have regular finite
element partitionings $\mathcal{T}_h^{i}$ of the subdomains
$\Omega_i$ into shape regular simplexes. These two meshes have corresponding 
trace meshes on the interface and for simplicity we assume that the meshes match across the interface so that
the trace meshes are equivalent and we may write
\begin{equation}\label{bmeshdef}
\mathcal{F}_h=\{ \ F \ : F=T\cap \Gamma, \ T\in \mathcal{T}_h^i, \; \text{$i=1$ or $2$}\ \}.
\end{equation}
We seek the approximation $u_h=(u_{1,h},u_{2,h})$ in the space
$V^h = V_{1}^h\times V_{2}^h$, where
\[
V^h_{i}=\left\{ v_{i}\in V_{i}:~v_{i}|_T\in \mathbb{P}_k(T), \,
\forall T \in \mathcal{T}_h \right\},\quad \text{for}\; k \ge 1,
\]
and for $\lambda_h$ 
in 
\[
\Lambda_h := \{q_h \in L_2(\Gamma): q_h\vert_F \in \mathbb{P}_l(F), \,
\forall F \in \mathcal{F}_h \},\quad \text{ for }\; l \ge 0 .
\]
If we now restrict the infimum of the supremum in (\ref{infsupref}) to our finite dimensional subspaces $V_h$ and
$\Lambda_h$ chosen such that the discrete spaces satisfy the
{\em inf-sup}\/ condition, we immediately obtain the standard Lagrange multiplier
domain decomposition method \cite{BrFo91}, with the Lagrangian given by
\[
\mathcal{L}(v_h,\mu_h) := \frac12a(v_h,v_h) +
\lga{\mu_h}{\jmp{v_h}} - \lom{f}{v_h} -\lga{g}{\avd{v_h}}.
\]
The augmented Lagrangian is obtained by adding a least squares penalty
on the constraint:
\[
\mathcal{L}(v_h,\mu_h) :=\frac12 a(v_h,v_h)+ 
\lga{\mu_h}{ \jmp{v_h}}+\frac{\gamma}{2}\|\jmp{v_h}\|^2_\Gamma
-\lom{f }{ v_h} -\lga{g}{\avd{v_h}}.
\]
This allows us to instead look for stationary points of the following
augmented Lagrangian
since we are working in discrete spaces and all pairings are $L^2$ scalar
products:

\begin{align*}
\mathcal{L}(v_h,\mu_h)   :={}&\frac12 a(v_h,v_h)+ 
\frac{1}{2 \gamma} \|\mu_h+\gamma  \jmp{v_h}\|^2_\Gamma \\
&
-\frac12 \lga{\gamma^{-1} \mu_h}{\mu_h} -\lom{f }{ v_h} -\lga{g}{\avd{v_h}}.
\end{align*}

The Euler--Lagrange equations characterising the saddle point of the
system (if it exists) takes the form: find $(u_h,\lambda_h) \in V_h
\times \Lambda_h$ such that
\begin{align}
a(u_h,v_h) + \gamma^{-1} \lga{\lambda_h + \gamma
  \jmp{u_h}}{\mu_h+\gamma \jmp{v_h}} - \lga{\gamma^{-1}
  \lambda_h}{\mu_h} = {}&\lom{f}{v_h} \nonumber\\
  & + \lga{g}{\avd{v_h}}
\end{align}
for all $v_h,\mu_h \in V_h \times \lambda_h$. Developing the second
term of the left hand side we see that this is equivalent to
\begin{align}\label{first_unstab}
a(u_h,v_h) + \lga{\lambda_h + \gamma
  \jmp{u_h}}{\jmp{v_h}} &= \lom{f}{v_h} + \lga{g}{\avd{v_h}} \\
\lga{ \jmp{u_h}}{\mu_h} & = 0.
\end{align}
This shows that the effect of the augmented Lagrangian compared to the
standard Lagrange multiplier method is simply the addition of a
penalty term on the constraint which, as mentioned above, gives us stronger control of
the constraints than would otherwise be possible. Indeed for the
standard Lagrange multiplier method only $\pi_\Lambda \jmp{u_h}$ is
controlled, where $\lga{\pi_\Lambda \jmp{u_h}}{\mu_h}
  =\lga{ \jmp{u_h}}{\mu_h} $ for all $\mu_h \in \Lambda_h$.
This formulation is however stable only for well balanced choices of $V_h$ and
$\Lambda_h$. In case the spaces do not satisfy the \emph{inf-sup}\/ condition
one may add a stabilizing term $j(\lambda_h,\lambda_h)$ satisfying
\[
\|h^{\frac12} (\lambda_h - \pi_\Gamma \lambda_h)\|_\Gamma \lesssim j(\lambda_h,\lambda_h)^{\frac12}
\]
where $\pi_\Gamma$ satisfies $\lga{\pi_\Gamma\mu_h}{v_h}
  =\lga{\mu_h}{v_h} $ for all $v_h$ in the trace mesh of
  $\mathcal{T}_1$ (or $\mathcal{T}_2$). For instance if $V_h$ is the
  space of piecewise affine, continuous functions in each subdomain
  and $\Lambda_h$ is
  the space of piecewise constant functions defined on the elements
  cut by $\Gamma$ we may choose
\[
j(\lambda_h,\lambda_h) := \sum_{F \in \mathcal{F}}
\|h\jmp{\lambda_h}\|_{\partial F\setminus \partial \Gamma}^2
\]
and follow the analysis of \cite{BuHa10a} to prove error estimates for
the formulation: find $(u_h,\lambda_h) \in V_h
\times \Lambda_h$ such that
\begin{align}\label{eq:mult_stab}
a(u_h,v_h) + \lga{\lambda_h + \gamma
  \jmp{u_h}}{\jmp{v_h}} &= \lom{f}{v_h} + \lga{g}{\avd{v_h}} \\
\lga{ \jmp{u_h}}{\mu_h} - j(\lambda_h,\mu_h)& = 0
\end{align}
for all $v_h,\mu_h \in V_h \times \lambda_h$.
On the other hand, if we formally replace $\lambda_h$ by
$-\av{\mu \partial_n u_h}$ and $\mu_h$ by $-\av{\mu \partial_n v_h}$
we obtain the formulation: find $u_h\in V_h$ such that
\begin{multline}\label{HHNitsche}
a(u_h,v_h) -\lga{\av{\mu \partial_n u_h}}{\jmp{v_h}} - \lga{ \jmp{u_h}}{\av{\mu \partial_n v_h}}  + \lga{ \gamma \jmp{u_h}}{\jmp{v_h}}\\= \lom{f}{v_h} + \lga{g}{\avd{v_h}}
\end{multline}
for all $v_h\in V_h$, and
we recognise Nitsche's formulation {from \cite{HaHa02}.}
\section{Debonding and adhesive contact\label{sec:adhesive}}
A robust discretization of the debonding problem was proposed in
\cite{HaHa04}. We will revisit their arguments in the context of
Lagrange multipliers as an augmented Lagrangian formulation
\label{sec:3}. The linear model problem in this case takes the form
\begin{align}\label{eq:debond_syst}
-\nabla \cdot \varepsilon \nabla u & = f \mbox{ in } \Omega_1 \cup
                                    \Omega_2\\
u&=0 \mbox{ on } \partial \Omega\\
\jmp{\varepsilon \partial_n u}&=0\mbox{ on } \Gamma\\
\jmp{u}&=-\compliance\av{\varepsilon \partial_n u} \mbox{ on } \Gamma. \label{eq:bond}
\end{align}
Here we note that by the continuity of the fluxes the formulation is
independent of the choice of the weights $w_i$, $i=1,2$. This time the
physical solution is discontinuous over the boundary and $u \in H^1(\Omega_1\cup\Omega_2)$.

The critical case is when $\compliance$ becomes large and a naive coupling
strategy leads to an ill-conditioned system or even locking on the
interface. We therefore follow \cite{HaHa04}, but contrary to the
discussion in {that} paper we herein use the augmented Lagrangian
formulation to arrive at the method. In the intermediate step we
obtain a robust Lagrange multiplier method for the debonding problem.
This time we start from the following Lagrangian, the saddle points
of which coincides with the solution of the debonding problem,
\[
\mathcal{L}(v,\mu) := \frac12 a(v,v) -
\left<\mu,\jmp{v} \right>_{-\frac12,\frac12,\Gamma}  - \frac12 \compliance \|\mu\|^2_\Gamma- \lom{f}{v}.
\]
The augmented Lagrangian is obtained adding a least squares term on
the constraint
\[
\mathcal{L}(v,\mu) := \frac12 a(v,v) -
\left<\mu,\jmp{v} \right>_{-\frac12,\frac12,\Gamma}  - \frac12 \compliance
\|\mu\|^2_\Gamma + \frac{1}{2\compliance} \|\jmp{v} +  \compliance \mu \|_\Gamma^2- \lom{f}{v}.
\]
Developing the square we see that the multiplier is eliminated and we
obtain
\begin{equation}\label{weq:start_lagrangian}
\mathcal{L}(v) := \frac12 a(v,v) + \frac{1}{2\compliance} \|\jmp{v} \|_\Gamma^2- \lom{f}{v}.
\end{equation}
Studying the corresponding optimality system leads to: find $u_h \in
V_h$ such that
\begin{equation}\label{weq:stiff_bonding}
a(u_h,v_h) + S \lga{\jmp{u_h}}{\jmp{v}} = \lom{f}{v},\quad \forall v_h
\in V_h.
\end{equation}
where $S=\compliance^{-1}$. As for the problem in Section \ref{sec:gen_bc} this formulation becomes ill-conditioned
for $\compliance$ small since the two terms of the left hand side will have sizes of
different orders of magnitude. A possible remedy is to replace $\compliance$ by
$\max(h,\compliance)$ in equation \eqref{weq:stiff_bonding}, but this results in
a nonconsistent perturbation of the system and reduced accuracy, in
the regime where $\compliance < h$.

If we instead modify the size of the least squares contribution in the
augmented Lagrangian we can moderate the strength of the imposition of
the constraint in a consistent manner. Introducing the parameter 
\[
S_h =(h_k/{\gamma_{\compliance}} + \compliance)^{-1},
\]
similarly as for the method in Section \ref{sec:gen_bc}, we may write
\[
\mathcal{L}(v,\mu) := \frac12 \|\varepsilon^{\frac12} \nabla v\|_\Omega^2 -
\left<\mu,\jmp{v} \right>_{-\frac12,\frac12,\Gamma}  - \frac12 \compliance
\|\mu\|^2_\Gamma + \frac{1}{2} S_h \|\jmp{v} +  \compliance \mu \|_\Gamma^2- \lom{f}{v}.
\]
Observe that the saddle point to this system is a weak solution to
\eqref{eq:debond_syst}--\eqref{eq:bond}, but the size of the weight in front of the
least squares term can never be larger than ${\gamma_{\compliance}}/h$, effectively
bounding the stiffness of the system. The corresponding optimality
system now reads:
\begin{multline}
a(u,v) - \lga{(1- \compliance S_h) \lambda}{\jmp{v}}-\lga{(1- \compliance S_h)
  \mu}{\jmp{u}}\\ - \lga{\compliance(1-S_h \compliance) \lambda}{\mu}
+\lga{S_h \jmp{u}}{\jmp{v}} = \lom{f}{v}.
\end{multline}
Formally replacing $\lambda$ by $\av{\varepsilon \partial_n u}$ and
$\mu$ by $\av{\varepsilon \partial_n v}$
and writing the form on the finite space $V_h$ yields: find $u_h \in
V_h$ such that
\begin{multline}
a(u_h,v_h) - \lga{(1- \compliance S_h) \av{\varepsilon \partial_n u_h}}{\jmp{v_h}}-\lga{(1- \compliance S_h)
  \av{\varepsilon \partial_n v_h}}{\jmp{u_h}}\\ - \lga{\compliance(1-S_h \compliance) \av{\varepsilon \partial_n u_h}}{\av{\varepsilon \partial_n v_h}}
\\ +\lga{S_h \jmp{u_h}}{\jmp{v_h}}  = \lom{f}{v_h}.\label{cohesive}
\end{multline}
We recognise the same Nitsche type method as proposed in \cite{HaHa04},
but this time with the weights chosen as in the previous section to
also be robust with respect to the contrast in the diffusivity.

\subsection{Adhesive contact}
In the previous linear model only the adhesive forces are accounted
for, which implies that penetration is possible. A more physically
realistic model excludes penetration by formulating the problem as a
variational inequality. Our last model problem concerns this nonlinear
model and we will combine the arguments developed above with those of
\cite{ChHi13,ChHiRe15}. Here for simplicity we assume that both
$\Omega_1$ and $\Omega_2$ intersects the boundary $\partial \Omega$.
\begin{align}
-\nabla \cdot \varepsilon \nabla u & = f \mbox{ in } \Omega_1 \cup
                                    \Omega_2 \label{eq:pde}\\
u&=0 \mbox{ on } \partial \Omega\\
\jmp{\varepsilon \partial_n u}&=  0\mbox{ on } \Gamma\label{eq:flux_jump}\\
\jmp{u}&  \leq 0\mbox{ on } \Gamma\label{eq:contact}\\
\compliance^{-1} \jmp{u}+\av{\varepsilon \partial_n u} & \leq 0\mbox{ on } \Gamma\\
\jmp{u}( \compliance^{-1}\jmp{u}+\av{\varepsilon \partial_n
  u}) & =0 \mbox{ on } \Gamma \label{eq:bondcont}
\end{align}
To cast this problem on the form of an augmented Lagrangian method we 
start out with the functional (\ref{weq:start_lagrangian}) with an additional constraint on $\jmp{v}$: 
\begin{equation}\label{eq:adh_cont_Lag}
\mathcal{L}(v,\mu) := \frac12 \|\varepsilon^{\frac12} \nabla
v\|_\Omega^2 -\left<\mu,\jmp{v}\right>_{-\frac12,\frac12,\Gamma}+ \frac{1}{2\compliance} \|\jmp{v} \|_\Gamma^2- \lom{f}{v}.
\end{equation}
The Euler--Lagrange equations  are to find $(u,\lambda)$ such that
\begin{equation}\label{one}
a(u,v)-\left<\lambda,\jmp{v}\right>_{-\frac12,\frac12,\Gamma} + \frac{1}{\compliance}\lga{ \jmp{u}}{\jmp{v}} = \lom{f}{v}\quad\forall v\in  H^1(\Omega_1\cup\Omega_2) ,
\end{equation}
\begin{equation}\label{wrongcons}
\left< \mu ,\jmp{u}\right>_{-\frac12,\frac12,\Gamma} = 0 \quad \forall \mu\in H^{-\frac12}(\Gamma) ,
\end{equation}
and we note from (\ref{one}) that, formally, the multiplier is given by
\begin{equation}\label{def}
\lambda = \av{\epsilon\partial_n u} + \compliance^{-1}\jmp{u}
\end{equation}
 but of course (\ref{wrongcons}) enforces $\jmp{u}=0$ weakly. In order to create a one-sided contact condition we
now consider the Kuhn--Tucker conditions (\ref{eq:contact})--(\ref{eq:bondcont}) as 
\begin{equation}
\lambda \leq 0, ~\jmp{u} \leq 0, ~\text{and}~  \lambda \jmp{u} = 0 \;\text{on $\Gamma$}. 
\end{equation}
These conditions can equivalently be formulated as
\begin{equation}\label{eq:adh_constraint}
\lambda = -{\gamma} [\jmp{u} -\gamma^{-1}\lambda]_+ ,
\end{equation}
where $[x]_+=\max (x,0)$, cf. \cite{ChHi13}. To introduce this
condition in the Lagrangian \eqref{eq:adh_cont_Lag} we note that, if
$\lambda,\mu \in L^2(\Gamma)$ (or if $\gamma^{-1}:H^{-\frac12}(\Gamma)\mapsto
H^{\frac12}(\Gamma)$, with suitable properties), in \eqref{one} we have
\[
-\lga{\lambda}{\jmp{v}} = -\lga{\lambda}{\jmp{v}-\gamma^{-1}\mu} -
\lga{\lambda}{\gamma^{-1}\mu}.
\]
Using now \eqref{eq:adh_constraint} in the first term of the right
hand side we have
\begin{equation}\label{right_const}
a(u,v)+\left<\gamma[\jmp{u} -\gamma^{-1}\lambda]_+,\jmp{v} -\gamma^{-1}\mu)\right>_\Gamma + \frac{1}{\compliance}\lga{ \jmp{u}}{\jmp{v}} = \lom{f}{v}.
\end{equation}
Moving over to discrete spaces, where the assumption on $\lambda_h$
makes sense, we write the corresponding Lagrangian on augmented form as
\[
\mathcal{L}_\text{a}(v_h,\mu_h) := a(v_h,v_h) +\frac{\gamma}{2}\Vert [\jmp{v_h}-\gamma^{-1}\mu_h]_+\Vert_\Gamma^2 - \frac{1}{2\gamma}\Vert \mu_h\Vert_\Gamma^2+\frac{1}{2 \compliance}\Vert \jmp{v_h}\Vert_\Gamma^2 - \lom{f}{v_h}
\]
leading to the problem of finding $(u_h,\lambda_h)\in V_h\times\Lambda_h$ such that
\begin{multline*}
a(u_h,v_h)+{\gamma}\left<[\jmp{u_h} -\gamma^{-1}\lambda]_+,\jmp{v_h}
  -\gamma^{-1}\mu_h)\right>_\Gamma -{\gamma}^{-1}\left< \lambda_h,
  \mu_h \right>_\Gamma\\
+ \frac{1}{\compliance}\left< \jmp{u_h}, \jmp{v_h} \right>_\Gamma =\lom{f}{v_h}
\end{multline*}
for all $(v_h,\mu_h)\in V_h\times\Lambda_h$.

Inserting now the definition (\ref{def}) of the multiplier and choosing $\mu = \av{\epsilon\partial_n v} + \compliance^{-1}\jmp{v}$ we seek $u_h\in V_h$ such that
\begin{align}\nonumber
{}& a(u_h,v_h) + \frac{1}{\compliance}\left< \jmp{u_h}, \jmp{v_h} \right>_\Gamma \\ \nonumber
{} &+{\gamma}\left<[(1-(\gamma \compliance)^{-1})\jmp{u_h} -\gamma^{-1}\av{\epsilon\partial_n u_h}]_+,(1-(\gamma \compliance)^{-1})\jmp{v} -\gamma^{-1}\av{\epsilon\partial_n v})\right>_\Gamma \\
{}& -\gamma^{-1}\left<\av{\epsilon\partial_n u_h} +
    \compliance^{-1}\jmp{u_h},\av{\epsilon\partial_n v_h} +
    \compliance^{-1}\jmp{v_h}\right>_\Gamma  = \lom{f}{v_h} 
 \quad\forall v_h\in V_h .\label{full_adhesion}
\end{align}

It is instructive to consider the two limiting cases of full contact and of no contact:
at contact, (\ref{full_adhesion}) gives
\begin{align*}
{}& a(u_h,v_h) -\left< \jmp{u_h},\av{\epsilon\partial_n v_h}\right>_\Gamma -\left< \jmp{v_h},\av{\epsilon\partial_n u_h}\right>_\Gamma \\
 {}&+ \left(\gamma-\frac{1}{\compliance}\right)\left<\jmp{u_h},\jmp{v_h}\right>_\Gamma= \lom{f}{v_h}\quad\forall v_h\in V_h.
\end{align*}
With the particular choice 
\[
\gamma := \frac{{\gamma_\compliance}}{h} + \frac{1}{\compliance},
\]
we obtain the following discrete problem: find $u_h\in V_h$ such that
\begin{equation*}
a(u_h,v_h)-\left<\av{\epsilon\partial_n u_h}\jmp{v_h}\right>-\left<\av{\epsilon\partial_n v_h}\jmp{u_h}\right> + \frac{{\gamma_\compliance}}{h}\left< \jmp{u_h}, \jmp{v_h} \right>_\Gamma = \lom{f}{v_h}
\end{equation*}
for all $v_h\in V_h$,
which is the standard Nitsche method {(\ref{HHNitsche}) for the adhesion free problem (with $\gamma_\compliance =\gamma_0$).} In the case of no contact observe
that by using $\gamma^{-1} = \gamma^{-1}-\compliance+\compliance$,
\begin{align*}
{}& \frac{1}{\compliance}\left< \jmp{u_h}, \jmp{v_h} \right>_\Gamma - \gamma^{-1}\left<\av{\epsilon\partial_n u_h} +
    \compliance^{-1}\jmp{u_h},\av{\epsilon\partial_n v_h} +
    \compliance^{-1}\jmp{v_h}\right>_\Gamma \\ 
{}\ & = -\left<\av{\epsilon\partial_n u_h},\jmp{v_h}+\compliance\av{\epsilon\partial_n v_h}\right>_\Gamma \\ {}&
 -\left<\jmp{u_h}+\compliance\av{\epsilon\partial_n u_h},\av{\epsilon\partial_n v_h}\right>_\Gamma -  \left<\compliance\av{\epsilon\partial_n u_h},\av{\epsilon\partial_n v_h}\right>_\Gamma \\ {}& 
(\compliance-\gamma^{-1}) \left<\av{\epsilon\partial_n u_h} +
    \compliance^{-1}\jmp{u_h},\av{\epsilon\partial_n v_h} +
    \compliance^{-1}\jmp{v_h}\right>_\Gamma .
\end{align*}
The finite element formulation then takes the form: we seek $u_h \in V_h$ such that
\begin{align*}
\lom{f}{v_h} = {}& a(u_h,v_h)-\left<\av{\epsilon\partial_n u_h},\jmp{v_h}+\compliance\av{\epsilon\partial_n v_h}\right>_\Gamma \\ {}&
 -\left<\jmp{u_h}+\compliance\av{\epsilon\partial_n u_h},\av{\epsilon\partial_n v_h}\right>_\Gamma -  \left<\compliance\av{\epsilon\partial_n u_h},\av{\epsilon\partial_n v_h}\right>_\Gamma\\
{}& + \frac{1}{\gamma_1}\left< \jmp{u_h}+\compliance\av{\epsilon\partial_n u_h}, \jmp{v_h} +\compliance\av{\epsilon\partial_n v_h}\right>_\Gamma , \quad\forall v_h\in V_h,
\end{align*}
where 
\[
\gamma_1 := \frac{\compliance^2}{\compliance-{\gamma^{-1}}} = \compliance+ \frac{h}{{\gamma_\compliance}},
\]
which coincides with the {form (\ref{cohesive}).}

\section{Stabilization for the extension to CutFEM\label{sec:cutmesh}}

Contrary to the methods discussed above, where the domains are meshed in the usual way, the CutFEM approach instead represents the boundary of a given domain 
on a background 
grid, for instance using a level set function. The background grid is 
then also used to represent the approximate solution of the governing partial differential 
equations. Consequently, CutFEM eases the burden of mesh generation by requiring only
a low-quality surface mesh representation of the
computational geometry.  Cutting the mesh will, however, result in boundary elements with very small
intersection with the physical domain. This may lead to a poorly
conditioned system matrix or failure of stability of the discrete
scheme. A remedy to this problem is to add a penalty term in the
cut element zone that extends the coercivity to the whole mesh domain,
i.e., in the $O(h)$ zone of the mesh domain (of each subdomain for interface problems) that does
not intersect the associated physical domain. This penalty term, termed \emph{ghost penalty}\/ due to it acting partly outside of the domain of interest, must
be carefully designed to add sufficient stability, while remaining
weakly consistent for smooth solutions. The basic methods described above are then applied to the cut meshes
and the only additional term is the ghost penalization.

To illustrate this idea, we consider the CutFEM method for the Poisson problem
\eqref{Nit1}.
We observe that by taking $v=u_h$ in the bilinear form $a(u_h,v)$, we have
the coercivity
\[
\|\nabla u_h\|_{L^2(\Omega)}^2 \leq a(u_h,u_h).
\]
However, to obtain coercivity of the form
$a_h(u_h,v)$ using this stability and the boundary penalty term, the
penalty parameter will depend on how the elements are cut, since, denoting the set of elements that are cut by $\Gamma$ by
\[
{\mathcal G}_h := \{K \in \mathcal{T}_h: K \cap \Gamma \neq \emptyset\}
\]
 we have, with
 \[
 a_h(u_h,v) :=a(u_h,v)-\lga{\partial_n u_h}{v}-\lga{\partial_n
   v_h}{u_h}+\lga{\gamma u_h}{v}
 \]
 that
\begin{align}
a_h(u_h,u_h) \ge {}& \|\nabla u_h\|^2_{L^2(\Omega)} + \|\gamma^{\frac12}  u_h\|^2_{L^2(\Gamma)}\nonumber \\ \label{Nitsche_coercivity}
& - 2 \sum_{K \in {\mathcal G}_h} \|\nabla u_h\|_{L^2(\Gamma \cap K)} \| u_h\|_{L^2(\Gamma \cap K)}.
\end{align}
Using the following well known
trace inequality:
under 
reasonable mesh assumptions
there exists a 
constant {$C_T$}, depending on $\Gamma$ but independent of the mesh, such that
\begin{equation}\label{indep1}
\| w\|_{L_2(\Gamma_K)}^2\leq C_T \left( h_K^{-1} \| w\|_{L_2(K)}^2 + h_K
\|\nabla w\|_{L_2(K)}^2\right), \quad \forall w\in H^1({K}) ,
\end{equation}
we have, 
\[
\|\nabla u_h\|_{L^2(\Gamma \cap K)}  \leq 
C_\Gamma \left(\frac{\vert K\cap\Gamma\vert}{\vert K \cap \Omega\vert}\right)^{\frac12} \|\nabla u_h\|_{L^2(K \cap \Omega)}, 
\]
where $\vert \cdot\vert$ denotes the measure of the indicated quantity. It follows that in principle we obtain coercivity by choosing 
\[
\gamma\vert_{K} >
2 {C_\gamma^2} \left(\frac{\vert K\cap\Gamma\vert}{\vert K \cap \Omega\vert}\right),
\] 
since by an arithmetic-geometric inequality, we have
\begin{align*}
a_h(u_h,u_h) \ge &{} \|\nabla u_h\|^2_{L^2(\Omega)} + \|\gamma^{\frac12} u_h\|^2_{L^2(\Gamma)}
- \frac12 \|\nabla u_h\|^2_{L^2(\Omega)} \\
& - \sum_{K \in {\mathcal G}_h} C_\Gamma^2 \left(\tfrac{\vert K\cap\Gamma\vert}{\vert K \cap \Omega\vert}\right)
\|u_h\|^2_{L^2(K \cap \Gamma)} \\
\ge &{}\frac12 \|\nabla
u_h\|^2_{L^2(\Omega)} + \|\left(\gamma - 2 C_\Gamma^2 \left(\tfrac{\vert K\cap\Gamma\vert}{\vert K \cap \Omega\vert}\right)\right)^{\frac12} u_h\|^2_{L^2(\Gamma)}  .
\end{align*}
Unfortunately this makes $\gamma$ strongly dependent on the cut, since
for $\vert K\cap\Gamma\vert = O(h_K)$, the volume measure $\vert K \cap \Omega\vert$ can be
arbitrarily small, resulting in problems both with conditioning and
accuracy. A solution
to this problem is to add a stabilizing term $g_h(u_h,v)$
to the form $a_h(\cdot,\cdot)$.
The role of this term is to extend
the coercivity from the physical domain $\Omega$ to the mesh domain
$\Omega_{\mathcal{T}}:= \Omega \cup \mathcal{G}_h$. In order to have this effect, the stabilization term should have the following properties.
\begin{enumerate}
\item It should give a bound on the energy norm on the mesh domain in the sense that
\begin{equation}\label{enhanced_coercivity}
c_G \|\nabla u_h\|^2_{\Omega_{\mathcal{T}}} \leq \|\nabla
u_h\|^2_{\Omega} + g_h(u_h,u_h),
\end{equation}
{where $c_G>0$ is bounded away from zero independent of the
mesh/boundary intersection for positive ghost penalty
stabilization parameter $\gamma_g$.}
\item  
For an interpolant of the extension of $u$, $i_h u := i_h \mathbb{E} u$ we must have the weak consistency
\[
g_h(i_h u,i_h u) \leq C h^k\|u\|_{H^{k+1}(\Omega)}.
\]
where the constant $C$ is independent of the mesh/boundary
intersection.
\end{enumerate}
One example of such a term is the ghost penalty stabilization
\begin{equation}\label{eq:ghostpenalty}
g_h(u_h,v) :=  \sum_{F \in \mathcal{F}_G} (\gamma_g h \jmp{\partial_{n_F} u_h},\jmp{\partial_{n_F} v})_{F},
\end{equation}
valid for piecewise affine approximation.
 Here, we introduced the set of element faces $\mathcal{F}_G$ associated with ${\mathcal G}_h$, defined as follows:
for each face $F \in \mathcal{F}_G$ there exists two simplices $K$ and $K'$ such that $F=K \cap K'$ and at least one of the two is a member of ${\mathcal G}_h$.
 This means in particular that the boundary faces of the mesh $\mathcal{T}_h$ are excluded from $\mathcal{F}_G$. We also used $\partial_{n_F}$ to denote the derivative in the direction of the normal to $F$.

Coercivity now follows from \eqref{Nitsche_coercivity}
and \eqref{enhanced_coercivity} as follows
\begin{align}\nonumber
a_h(u_h,u_h) \ge &{} \|\nabla u_h\|^2_{L^2(\Omega)} + \gamma \|h^{-\frac12} u_h\|^2_{L^2(\Gamma)} \\ \nonumber
& - 2 C_T \|\nabla u_h\|_{L^2(\mathcal{G}_h)} \|h^{-\frac12} u_h\|_{L^2(\Gamma)} +
g_h(u_h,u_h) \\ \nonumber
 \ge &{} c_G \|\nabla u_h\|^2_{\Omega_{\mathcal{T}}}   + \gamma
\|h^{-\frac12} u_h\|^2_{L^2(\Gamma)}- 2 C_T \|\nabla
u_h\|_{L^2(\Omega_{\mathcal{T}})} \|h^{-\frac12} u_h\|_{L^2(\Gamma)}  \\ 
\ge &{} \frac{c_G}{2}  \|\nabla u_h\|^2_{\Omega_{\mathcal{T}}} + (\gamma
- 2 C_T^2 c_G^{-1}) \|h^{-\frac12} u_h\|^2_{L^2(\Gamma)}.\label{full_coercivity}
\end{align}
Here $C_T$ is the constant of the trace inequality \eqref{indep1} and $c_G$
is the coercivity constant of the stability estimate \eqref{enhanced_coercivity}.
We conclude by choosing $\gamma_0 > 2 C_T^2 c_G^{-1}$, where the
lower bound is independent of the mesh/boundary intersection, but not
of the penalty parameter $\gamma_g$ in $g_h(\cdot,\cdot)$. Error
estimates now follow in a similar fashion as for the standard
Nitsche's method, using \eqref{full_coercivity} and
the consistency of the penalty term.
One may also show that the conditioning of the
system matrix is bounded independently of the mesh/boundary
intersection. For further details see \cite{BuHa12}.

Extension to the other model problems is straightforward. We consider 
the problem of interface coupling using \eqref{HHNitsche}. In the original
paper on cut finite elements \cite{HaHa02} the method for meshed subdomains was
carried over to the cut element case using piecewise affine elements and weights
\begin{equation}\label{geomweight}
w_1 = K\cap \Omega_2/|K|,\quad w_2 = K\cap \Omega_1/|K| .
\end{equation}
However, for problems with large contrast $\varepsilon_{max}/\varepsilon_{min}$
it is known that thise choice is not stable for
arbitrary cuts if the mesh size is not small enough to resolve the
contrast. Indeed too large
contrast can lead to a phenomenon reminiscent of locking {for}
unfortunate cuts (i.e. if no $H^1$-conforming subspace with
approximation exists). In case robustness is necessary 
we instead choose the weights to be $w_1 =
\epsilon_2/(\epsilon_1+\epsilon_2)$ and $w_2 =
\epsilon_1/(\epsilon_1+\epsilon_2)$ and $\gamma^{-1} = \gamma_0 h^{-1}
\omega(\epsilon)$ with $\omega(\epsilon):=2 (\epsilon_1 \epsilon_2)/(\epsilon_1+\epsilon_2)$ we identify this method as the
Nitsche method discussed in \cite{BuZU06, BuZu10} which was shown to be
stable on unfitted meshes provided a ghost penalty term
is added. This additional stabilization term should here have the
properties analogous to those for the cut fictitious domain method discussed above, and
a typical example is the modification of \eqref{eq:ghostpenalty} now acting across the faces of the cut elements on the interface.
The resulting method takes the form
\begin{multline}\label{BZNitsche}
a(u_h,v_h) -\lga{\av{\mu \partial_n u_h}}{\jmp{v_h}} - \lga{ \jmp{u_h}}{\av{\mu \partial_n v_h}}  \\
+ \lga{ \gamma \jmp{u_h}}{\jmp{v_h}}+g_h(u_h,v_h)\\= \lom{f}{v_h} + \lga{g}{\avd{v_h}}.
\end{multline}
In particular it was proven in \cite{BGSS16} that for piecewise affine
elements and smooth enough
$\Gamma$ there holds
\[
\sum_{i=1}^2 \|\varepsilon_i \nabla (u_i - u_{i,h})\|_{\Omega_i} \leq C h \|f\|_{\Omega} 
\]
where the constant $C$ is independent of $\varepsilon$ and $h$. This
result is possible to obtain thanks to the fact that the weights
shift the interface term to the side where $\varepsilon_i$ is the
smallest. 

Observe now the resemblance between the formulation \eqref{BZNitsche}
and \eqref{eq:mult_stab}. The latter formulation however is not in
general robust
for large contrast.
Indeed regardless of the contrast the constraint will be satisfied
equally strongly and we know from the experience of \cite{BZ12}
that the trick to obtaining robustness is to relax the control
obtained by the multiplier by redefining $\lambda_h$. Instead of
identifying $\lambda_h = \av{\varepsilon \partial_n u_h}$ we use that
$\av{\varepsilon \partial_n u_h} = \omega(\varepsilon) \{\partial_n
  u_h\}$ where $\{\cdot\}$ denotes the standard
  arithmetic average. If we then instead identify $\lambda_h= \{\partial_n
  u_h\}$ and introduce the factor $ \omega(\varepsilon)$ in the
  formulation \eqref{eq:mult_stab} we obtain 
\begin{align}\label{eq:mult_stab4}
a(u_h,v_h) + \lga{\omega(\varepsilon) \lambda_h + \gamma
  \jmp{u_h}}{\jmp{v_h}} &= \lom{f}{v_h} + \lga{g}{\avd{v_h}} \\
\lga{ \jmp{u_h}}{\omega(\varepsilon) \mu_h} -j(\lambda_h,\mu_h)& = 0
\end{align}
where once again $\gamma = \gamma_0 h^{-1}
\omega(\epsilon)$ and the stabilization operator $j(\cdot,\cdot)$ also
must scale as $\omega(\varepsilon)$. This formulation will relax the
jump in a similar fashion as \eqref{BZNitsche}, but control of
$\lambda_h$ is sacrificed as $\omega(\varepsilon)$ becomes small.

\section{A numerical example\label{sec:example}}

We give an example of how the method works in the case of adhesion, with and without contact.
To exemplify how the different aspect of adhesion come into play we consider 
a domain $(0,1)\times(0,1)$ with $u=0$ at $x=0$ and at $y=0$ and with $\partial_n u=0$ on other boundaries. The domain is cut by a half circle with radius $r=0.74$; $\varepsilon =2$ on the domain $\Omega_1$ containing the origin and $\varepsilon =  1/2$ on $\Omega_2$. The right-hand side is given by
\[
f= \left\{\begin{array}{c} 1\;\; \text{if $y\leq 1/2$,}\\[3mm]-7/2\;\; \text{if $y > 1/2$.}\end{array}\right.
\]
We set $\gamma_0=100$, used no jump stabilization, and set $\omega(\varepsilon) = 1$ {but used geometric averages of the type (\ref{geomweight}).}

\begin{figure}[hbt]
\begin{center}
\includegraphics[width=3in]{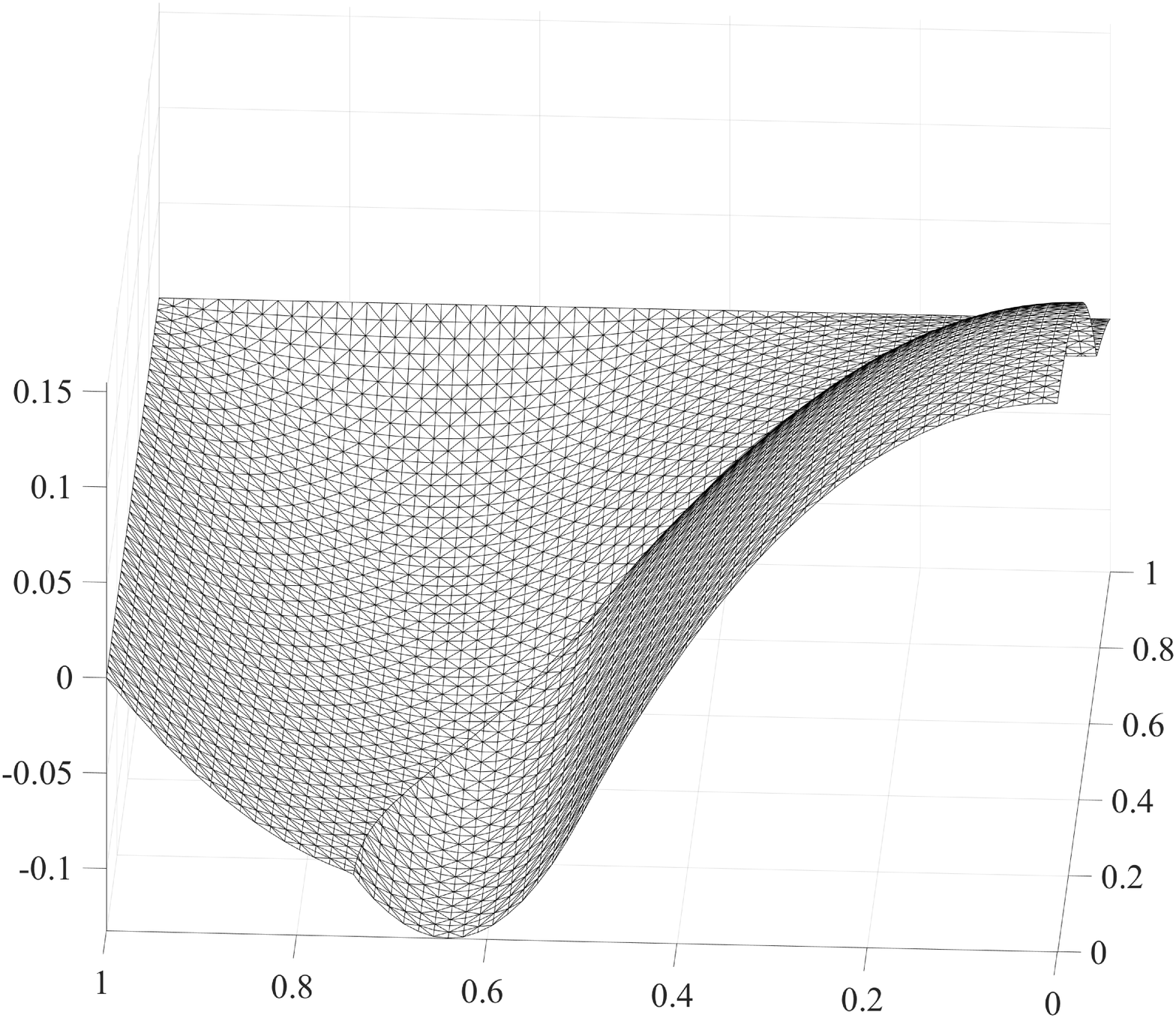}
\end{center}
\caption{Continuity enforced by Nitsche's method.\label{fig:cont}}
\end{figure}

\begin{figure}[hbt]
\begin{center}
\includegraphics[width=3in]{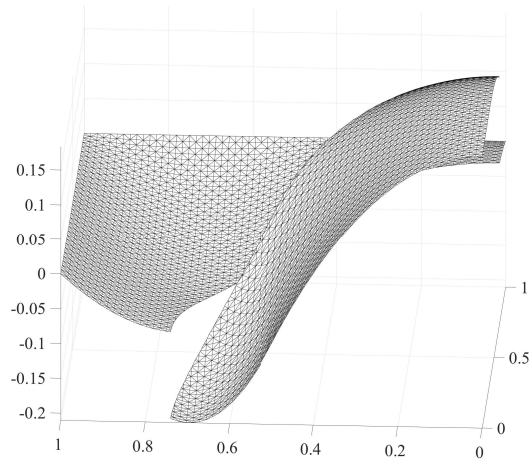}
\end{center}
\caption{A cohesive interface law enforced by Nitsche's method.\label{fig:discont}}
\end{figure}

\begin{figure}[hbt]
\begin{center}
\includegraphics[width=3in]{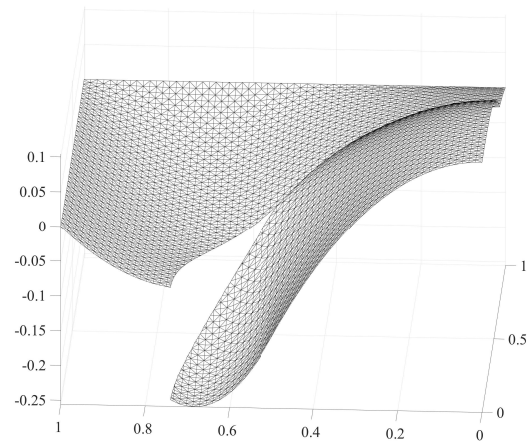}
\end{center}
\caption{Cohesive interface combined with a contact condition by Nitsche's method.\label{fig:contact}}
\end{figure}

 In Fig.\ref{fig:cont} we show the solution using the standard Nitsche method (\ref{HHNitsche}),
in Fig. \ref{fig:discont} we show the solution with a cohesive interface with $\compliance=1/2$ using (\ref{cohesive}), and in Fig. \ref{fig:contact} we show the
solution using a one sided contact condition as in (\ref{full_adhesion}).

\section*{Acknowledgements}
The contribution of the first author was supported in part by the EPSRC grants EP/J002313/2 and
EP/P01576X/1, the contribution of the
second author was supported in part by the Swedish Foundation for Strategic Research Grant No. AM13-0029 and the Swedish Research Council Grant No. 2011-4992.
%
%

\newpage

\end{document}